\documentclass[letterpaper,10pt]{article}
\usepackage[utf8]{inputenc}
\usepackage{amsmath, amsthm, amsfonts, amssymb}
\usepackage{mathabx}
\usepackage{todonotes}
\usepackage{stmaryrd}
\usepackage{graphicx}
\usepackage{tikz}
\usetikzlibrary{calc}
\usetikzlibrary{shapes}
\usepackage{pgfplots}
\pgfplotsset{compat=newest}
\usepackage{microtype}
\title{Classification of Quadratic Packing Polynomials on Sectors of $\R^2$}
\author{Madeline Brandt \and K{\aa}re Gjaldb{\ae}k}

\newtheorem{thm}{Theorem}
\newtheorem{lem}[thm]{Lemma}
\newtheorem{prop}[thm]{Proposition}
\newtheorem*{fact*}{Fact}

\theoremstyle{definition}
\newtheorem{rem}[thm]{Remark}
\newtheorem{dfn}[thm]{Definition}
\newtheorem{example}[thm]{Example}
\newtheorem{question}{Question}

\newcommand{\N}{\mathbb{N}} 
\newcommand{\Z}{\mathbb{Z}} 
\newcommand{\Q}{\mathbb{Q}} 
\newcommand{\R}{\mathbb{R}} 
\newcommand{\forwhich}{~:~} 
\newcommand{\il}{\textrm} 
\newcommand{\size}[1]{\#{#1}} 
\newcommand{\yif}[1]{\bar{y}_{#1}}	

\newenvironment{case}[1]{\par\noindent\textit{Case: #1}\newline\noindent}{}

\begin{document}

\maketitle

\begin{abstract}
    We study quadratic polynomials giving bijections
    from the integer lattice points of sectors of $\R^2$ onto $\N_0$,
    called packing polynomials.
    We determine all quadratic packing polynomials on rational sectors.
    This generalizes results of Stanton \cite{stan14}, Nathanson \cite{nath14},
    and Fueter and P{\'o}lya \cite{fuetPol23}.
\end{abstract}



\section*{Introduction}
    In the 1870s, Cantor introduces the notion of comparing the sizes of sets using
    bijective correspondences \cite{cantor1878beitrag}.
    He observes that the polynomial
    \begin{equation*}
        f(x,y) = x + \frac{(x+y-1)(x+y-2)}{2}
    \end{equation*}
     has the \emph{remarkable property}%
\footnote{\emph{bemerkenswerthe Eigenshaft}, in the words of the author.}
    that it represents each positive integer exactly once when $x$ and $y$
    run over all positive integers.
    In \cite{fuetPol23}, Fueter and P{\'o}lya prove that two derived
    \emph{Cantor polynomials}%
\footnote{
    Adjusted to include $0$, one is obtained from the other by swapping coordinates.
    Lew and Rosenberg \cite{lewRos78} use the term \emph{Cantor polynomials}.
    Other authors (e.g. Smory{\'n}ski \cite{smor1991lnt}) use the term
    \emph{Cantor pairing functions}. Fueter and P{\'o}lya leave them unnamed.
}
  \begin{equation}\label{eq:cantor}
        F(x,y) = \frac{1}{2}(x+y)(x+y+1) + x,\ \ \ \  G(x,y) = \frac{1}{2}(x+y)(x+y+1) + y
  \end{equation}
    are the only \emph{quadratic} polynomials which bijectively map $\N_0 \times \N_0$ to $\N_0$.
    The authors also formulate the \emph{Fueter-P{\'o}lya Conjecture} which states that
    the Cantor polynomials are the only polynomials of any degree which admit such a bijection.
    
    Lew and Rosenberg developed the theory further and studied the Fueter-P{\'o}lya Conjecture in the 1970s \cite{lewRos78,lewRos78p2}.
    They prove that polynomial bijections of degrees $3$ and $4$ are impossible, but the general conjecture remains open.
    They introduce the term \emph{packing function} for a function
    which bijectively maps a discrete set to the non-negative integers.
    For some preliminary results, they study the lattice points of 
    \emph{sectors} of $\mathbb{R}^2$, which are sets $\{(x,y) \in \mathbb{R}^2\ |\ y\leq \alpha x,\ x\geq 0,\ y\geq 0\}$. Allowing for $\alpha = \infty$, we may think of the first quadrant as a special case of this.
    
    Nathanson studies quadratic packing polynomials on the lattice points of sectors of $\mathbb{R}^2$ in \cite{nath14}.
    He finds two packing polynomials for each sector with $\alpha = \frac{1}{m}$, $m \in \N$.
    He shows that
    \begin{align*}
        F_{\frac{1}{m}}(x,y) &= \frac{1}{2}(x-(m-1)y)(x-(m-1)y - 1) + x + (1-(m-1))y,\\
        G_{\frac{1}{m}}(x,y) &= \frac{1}{2}(x-(m-1)y)(x-(m-1)y + 1) + x + (-1-(m-1))y
    \end{align*}
    are the only two quadratic packing polynomials on sectors of this type.
    For each sector defined by $\alpha = n$, $n\in \N$,
    he finds two packing polynomials
    \begin{equation}\label{eq:nathintpols}
        F_{n}(x,y) = \frac{n}{2}x(x-1) + x + y,\ \ \text{and} \ \
        G_{n}(x,y) = \frac{n}{2}x(x+1) + x - y,
    \end{equation}
    but does not rule out the possibility of more.
    Nathanson's article ends with a list of open problems:
    
    \begin{question}\label{question1}
        Are there packing polynomials on sectors with irrational $\alpha$?
    \end{question}
    \begin{question}\label{question2}
        Are there packing polynomials of higher degree on sectors?
    \end{question}
    \begin{question}\label{question3}
        Classify quadratic packing polynomials on sectors with $\alpha \in \Q$.
    \end{question}

    Question~\ref{question1} is settled in
    the quadratic case by Gjaldb{\ae}k \cite[Corollary 6]{gjald20},
    and independently by Sury and Vsemirnov \cite{sury2022packing},
    while Question~\ref{question2} remains open.
    Stanton \cite{stan14} studies Question~\ref{question3}.
    She provides a necessary form for the homogeneous quadratic part
    of the polynomial and uses it to classify all quadratic packing polynomials
    for integral $\alpha$.
    Shortly thereafter, Brandt \cite{brandt14} finds a method to tackle the general rational case.
    The article is incomplete and remained in preprint, but the method is fruitful and
    the conclusions correct, as we shall see. 
    In this article, we provide a complete classification of quadratic packing polynomials on sectors with rational slope (Theorem \ref{thm:main}), thus answering Question~\ref{question3} and completing the classification of quadratic packing polynomials on sectors.
    As special cases we recover the integral case (Stanton \cite{stan14}) and 
    the reciprocals (Nathanson \cite{nath14}).
    The original scenario (Fueter-P{\'o}lya \cite{fuetPol23}),
    while not included as a case in Theorem \ref{thm:main}, follows from it.
    
\section{Quadratic Packing Polynomials on Sectors}
	
	In this section we define sectors and quadratic packing polynomials. Then, we state Stanton's necessary form for a quadratic packing polynomial, and provide an independent proof of this result.
	
	For two linearly independent vectors $\omega_1, \omega_2 \in \R^2$, define
	the sector $S(\omega_1, \omega_2)$ as the conical hull of $\omega_1$ and $\omega_2$.
	That is
	\[
	    S(\omega_1, \omega_2) := \{ t_1 \omega_1 + t_2 \omega_2 \forwhich t_1, t_2 \geq 0 \}.
	\]
	Alternatively, $S(\omega_1, \omega_2)$ is the convex hull of the rays
	$t \omega_1$, $t\omega_2$, $t \geq 0$.
	We will restrict our attention to sectors where one of the rays is the positive $x$-axis. As long as one of $\omega_1$ or $\omega_2$ is rational, this does not present any loss of generality (see Remark \ref{rem:omega}).
	We denote for $\alpha > 0$ the sector $S(\alpha) \subset \R^2$ as
	\[
		S(\alpha) := S\big((1,0), (1,\alpha)\big) = \{(x,y) \in \R^2 \forwhich 0 \leq y \leq \alpha x\}
	\]
	and the integer lattice points in $S(\alpha)$ as
	\[
		I(\alpha) := S(\alpha) \cap \Z^2.
	\]
	$S(\infty)$ and $I(\infty)$ will denote the first quadrant and its lattice points, respectively. To be precise
    \[
        S(\infty) := S((1,0),(0,1)), \quad I(\infty) := S(\infty) \cap \Z^2.
    \]
	We will refer to sectors as integral, rational, irrational according to $\alpha$ being
	integral, rational, irrational. The sector $S(\infty)$ is considered integral.

	\begin{dfn}
	    Let $I \subset \mathbb{R}^2$ be an enumerable set.
		A function $f:\R^2\to\R$ is called a \emph{packing function} on $I$
		if it maps $I$ bijectively onto $\N_0$.
		If $f$ is a polynomial it is called a \emph{packing polynomial}.
	\end{dfn}
	
	\begin{rem}
	    In Lew and Rosenberg's definition the set $I$ is always a subset of $\Z^m$
	    and Nathanson uses this definition as well.
	    Our broader definition is practical for the treatment we present here.
	    If we speak of a packing polynomial on $S(\alpha)$ the enumerable set implied
	    will be $I(\alpha)$.
	\end{rem}
	
	In this paper we will classify all quadratic packing polynomials (QPPs) 
    on the lattice points of sectors $S(\alpha)$.
    We will determine which values $\alpha$ allow for QPPs and provide
    formulas for each.
    
    \begin{example} 
        \label{ex:43}
        Let $\alpha = 4/3$. Then the polynomial
        $
            2 x^2-2 x y+\frac{y^2}{2}+\frac{y}{2}
        $
        is a quadratic packing polynomial on $S(4/3)$ by Theorem \ref{thm:main}, see Figure \ref{fig:43}.
        \begin{figure}[htbp]
            \centering
            \includegraphics{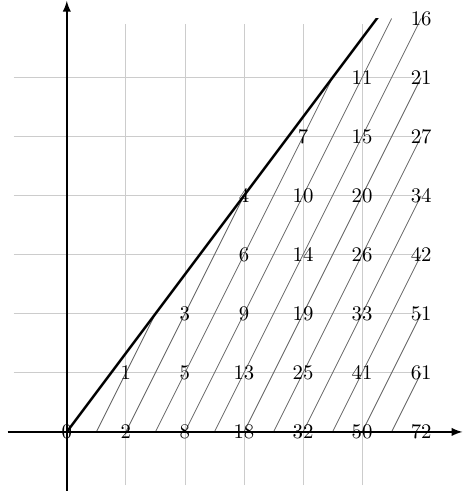}
            \caption{The quadratic packing polynomial on $S(4/3)$ given in Example \ref{ex:43}.}
            \label{fig:43}
        \end{figure}
   \end{example}
    
    Our starting point is the following result (see \cite{gjald20}, Equation~1 and Theorem~5).
	\begin{thm}[Gjaldb{\ae}k \cite{gjald20}]\label{thm:alphaForm}
		If $P(x,y)$ is a QPP on $I(\alpha)$, then
		\[
			P(x,y) = \frac{A}{2} x(x-1) + B xy + \frac{C}{2} y(y-1) + Dx + Ey + F
		\]
		with $A, B, C, D, E, F \in \Z$.
		Furthermore, we must have $A > 0$, $B^2 = AC$ and
		\[
			\alpha = \frac{A}{1-B},
		\]
		so in particular, $\alpha$ is rational.
	\end{thm}
	
	Since $\alpha$ must be rational in order for a packing polynomial on $S(\alpha)$ to exist, we focus on the case where $\alpha$ is rational in this article.
	
    Stanton gives the coefficients of
    the homogeneous quadratic part of $P(x,y)$
    in terms of the slope of the sector and a necessary condition on the slope \cite{stan14}.
	
	\begin{thm}[Stanton \cite{stan14}]\label{thm:stanNecForm}
		Let $\alpha = \frac{n}{m}$ be rational with $\gcd(m,n) = 1$.
		If $P(x,y)$ is a QPP on $I\left(\frac{n}{m}\right)$, then $n \mid (m-1)^2$ and
		\[
			P(x,y) = \frac{n}{2}\left(x - \frac{m-1}{n} y\right)^2 
				+ \left(D - \frac{n}{2}\right) x+ \left(E - \frac{(m-1)^2/n}{2}\right) y + F
		\]
		with $D, E, F \in \Z$.
	\end{thm}
    This was proved using the Lindemann-Weierstra{\ss} Theorem.
	Using Theorem~\ref{thm:alphaForm},
	we are able to present a short alternative proof.
	\begin{proof}
		By Theorem \ref{thm:alphaForm},
		\[
			\alpha = \frac{n}{m} = \frac{A}{1-B}, \il{ so } n(1-B) = mA.
		\]
		Since $m$ and $n$ are coprime, we have $n \mid A$.
		Since $AC = B^2$ by Theorem \ref{thm:alphaForm},
		\[
			nAC = nB^2 = n(1-(A/n)m)^2 = n - A(2  m - (A/n)  m^2),
		\]
		so that $n = A(nC + 2  m - (A/n)  m^2)$
    	so $A \mid n$. Since $A > 0$, we conclude that $A = n$.
		It follows that $B = 1 - (A/n)m = 1 - m$ and
		$C = B^2/A = (1-m)^2/n$. That
		$C$ is an integer provides the condition $n \mid (m-1)^2$.
	\end{proof}
	
	\begin{rem}
	    The case of $I(\infty)$ corresponds to $n=1$, $m=0$.
        We will elaborate on this later.
	\end{rem}

    \begin{example}\label{ex:3over2}
        Since $3 \nmid (2-1)^2$, there are no QPPs on the sector $I\left(\frac{3}{2}\right)$.
    \end{example}
	
\section{$k$-stairs and Sector Transformation}

    In this section, we introduce two tools that we will use to prove our main result.
    The first is the notion of \emph{staircases},
    which help us describe the structure of quadratic packing polynomials.
    Next, we describe how to transform sectors so that we can more readily leverage the structure provided by staircases.

\begin{figure}[ht]
	\centering
    \includegraphics{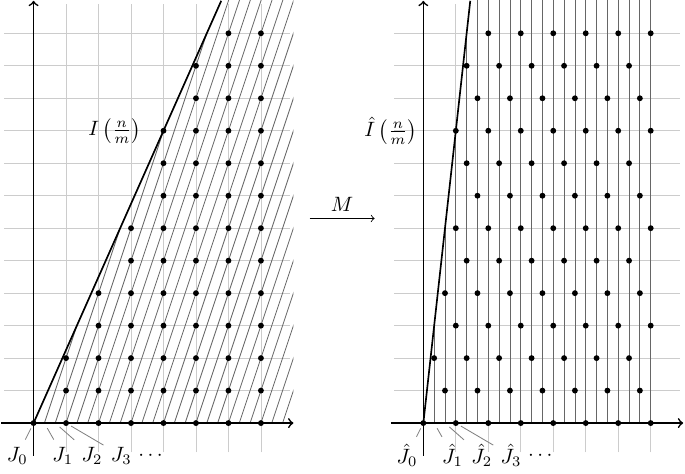}
    \caption{Transforming the sector $S\left(\frac{n}{m}\right)$.
        Example showing $\frac{n}{m} = \frac{9}{4}$.}
	\label{fig:trfSector}
\end{figure}

	The points in $I\left(\frac{n}{m}\right)$ can be
	subdivided into disjoint sets of points that fall on the same line with slope $\frac{n}{m-1}$
	(or vertical lines in the case $m=1$), see the left side of Figure \ref{fig:trfSector}.
	Let $l = \gcd(m-1, n)$.
	For an integer $i$ we call
	\[
		J_i := \left\{
			(x,y) \in I\left(\frac{n}{m}\right) \forwhich
			(m-1)y = n x - l i
		\right\}
	\]
	the \emph{$i$th staircase}, noting that $I\left(\frac{n}{m}\right)$ is the disjoint union
	\[
		I\left(\frac{n}{m}\right) = \bigsqcup_{i=0}^{\infty} J_i.
	\]
	We refer to the points as \emph{steps} on the staircase
	with the first step being the one with the minimal $y$-coordinate.
	For consecutive steps $(x_1, y_1)$ and $(x_2, y_2)$ on a staircase, we have
	$(x_2, y_2) = (x_1 + (m-1)/l, y_1 + n/l)$.
	The difference of $P(x,y)$ evaluated on consecutive steps, regardless of the staircase, is
	\begin{align*}
		&P(x + (m-1)/l, y + n/l) - P(x,y)\\
			&\hspace{2cm}=\left(D-\frac{n}{2}\right) (m-1)/l + \left(E - \frac{(m-1)^2/n}{2}\right) n/l\\
			&\hspace{2cm}:= k.
	\end{align*}

    We will now transform sectors so that the staircases have a nicer structure.
	For an invertible $2\times 2$ matrix $M$ we study the transformed sector $M(S(\alpha))$ instead of $S(\alpha)$.
	If $f(x,y)$ is a packing function on $I(\alpha)$, then
	$f(M^{-1}(x,y))$ is a packing function on $M(I(\alpha))$. Likewise,
	any packing function on the transformed lattice corresponds to a packing function
	on the original.%
\footnote{
    Nathanson \cite{nath14} used transformations of the form
    $\begin{pmatrix}
		1 & m \\
		0 & 1
	\end{pmatrix}$
	to obtain equivalences between the sector $S(\infty)$
	and sectors $S\left(\frac{1}{m}\right)$.
}
	Consider the transformation
	\[
		M = \begin{pmatrix}
				1 & -\frac{m-1}{n} \\
				0 & 1
		    \end{pmatrix},\quad
		M^{-1} = \begin{pmatrix}
				1 & \frac{m-1}{n} \\
				0 & 1
		    \end{pmatrix}.
	\]
	This transformation skews the sector and the lattice points with it, turning all staircases vertical, see Figure \ref{fig:trfSector}.
	Note that in the case of integral sectors the staircases are already vertical
	and the transformations are the identity.
	The sector $S(\infty)$ is sent to $S(1)$, so there is really no need to treat $S(\infty)$
	as a special case.
    We will use the following notation:
	\begin{align*}
    	\hat{S}\left(\frac{n}{m}\right) &:= M\left(S\left(\frac{n}{m}\right)\right),\quad
		\hat{I}\left(\frac{n}{m}\right) := M\left(I\left(\frac{n}{m}\right)\right),\quad
		\hat{J}_i := M\left(J_i\right).
	\end{align*}
    \begin{rem}\label{rem:trfSectorStructure}
    	Note that $\hat{S}\left(\frac{n}{m}\right) = S(n)$
    	and $\hat{I}\left(\frac{n}{m}\right) = M\left(\Z^2\right) \cap S(n) = \bigsqcup_{i=0}^{\infty} \hat{J}_i$.
        If $(a,b) \in \hat{J}$, then
        $\left(a + \frac{m-1}{n} b, b\right)$
        is an integer coordinate point on the line
        $((m-1)/l) y = (n/l)x - i$,
        so
        \[
            a = \frac{i}{n/l},\quad\il{and}\quad b \equiv - \frac{1}{(m-1)/l} i \pmod {n/l}
        \]
    	Also note that the number of points on $\hat{J}_i$ is given by
    	\[
    	    \size{\hat{J}_i} = (l^2/n) i + \llbracket n/l \mid i \rrbracket,
    	\]
    	where
    	\[
    	    \llbracket n/l \mid i \rrbracket = \left\{
    	        \begin{array}{ll}
    	            0 & \il{if $n/l$ does not divide $i$},\\
    	            1 & \il{if $n/l$ divides $i$}.
    	        \end{array}
    	        \right.
    	\]
    \end{rem}

    \begin{rem}\label{rem:trfSectorInfty}
        The sectors $S(\infty)$ and $S(1)$ are equivalent by
        the transformations
        \[
            M = \begin{pmatrix}
				1 & 1 \\
				0 & 1
		    \end{pmatrix},\quad
		      M^{-1} = \begin{pmatrix}
				1 & -1 \\
				0 & 1
		    \end{pmatrix}.
        \]
        Remembering that the rational sector $S\left(\frac{n}{m}\right)$
        is shorthand notation for the equivalent
        $S\left( (1,0), (m,n) \right)$,
        and $S(\infty)$ is short for $S\left( (1,0), (0,1) \right)$,
        it makes sense for convenience to consider $S(\infty)$ the rational
        sector with $n = 1$ and $m = 0$.
    \end{rem}
    
	The vertical staircases simplify computations. We set
	\begin{equation}\label{eq:initNecForm}
		\hat{P}(x,y) := P\left(M^{-1}(x,y)\right) = \frac{n}{2} x (x-1) + D x + \frac{k}{n/l} y + F.
	\end{equation}
	$\hat{P}(x,y)$ is a QPP on $\hat{I}\left(\frac{n}{m}\right)$
	if and only if $P(x,y)$ is a QPP on $I\left(\frac{n}{m}\right)$.
	
\begin{example}\label{ex:throughEx1}
    The polynomial
    \[
        P(x,y) = 6\left(x - \frac{1}{2}y\right)
                    \left(x - \frac{1}{2}y - \frac{3}{2}\right)
                    + x + y + 2
    \]
    is a QPP with $k = 3$ on $I\left(\frac{12}{7}\right)$,
    according to Theorem \ref{thm:main}.
    The corresponding QPP on $\hat{I}\left(\frac{12}{7}\right)$
    is
    \[
        \hat{P}(x,y) = 6x \left(x-\frac{3}{2}\right) + x + \frac{3}{2} y + 2.
    \]
    See Figure \ref{fig:12o7kp3Ex1}.

    \begin{figure}
        \centering
        \includegraphics[scale=0.75]{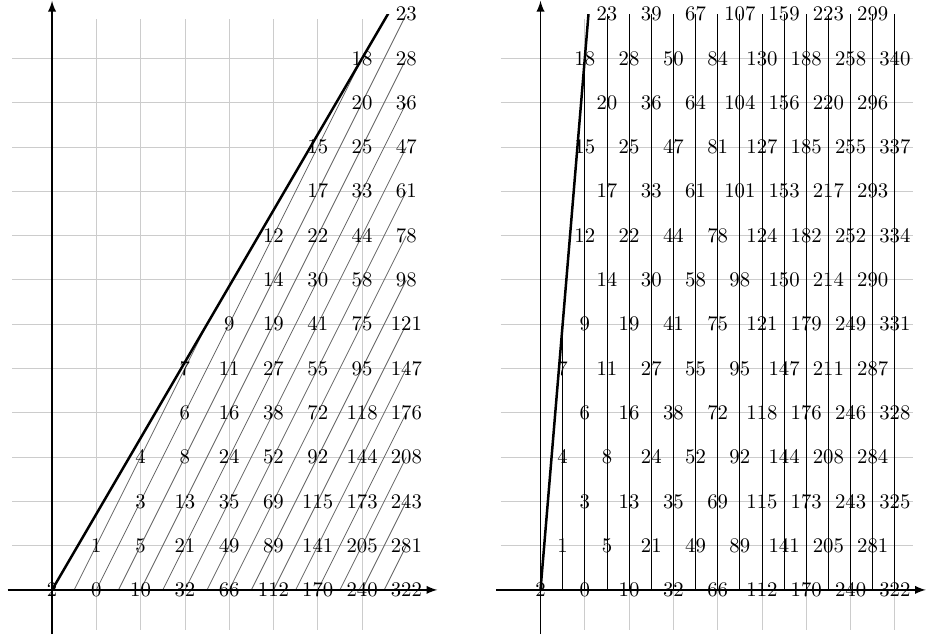}
        \caption{A QPP on the sector $I\left(\frac{12}{7}\right)$
        and the transformed sector $\hat{I}\left(\frac{12}{7}\right)$.}
        \label{fig:12o7kp3Ex1}
    \end{figure}
\end{example}

\section{Necessary Conditions}

In this section, we prove a series of lemmas that describe the structure of QPPs on transformed sectors. These lead to Proposition \ref{prop:necCond}, which gives the necessary form of $\hat{P}(x,y)$.

    Recall that $k$ denotes the difference $P(x_2,y_2) - P(x_1,y_1)$ of $P$ on consecutive stairs.
    Initially, we will assume that $k>0$. 
    Later, in the proof of Theorem \ref{thm:main}, we will see that the case $k<0$ easily follows.
	Let $\yif{i}$
    denote the $y$-coordinate of the first step on $J_i$.
    Note that in the case of integral sectors much of the discussion below
    could be greatly simplified, since the sector transformation is trivial
    and $\yif{i} = 0$ for all $i$.
	
	\begin{lem}\label{lem:nextCong}
	    Let $\hat{P}(x,y)$ be a QPP on $\hat{I}\left(\frac{n}{m}\right)$ with $k > 0$.
		Then
		\[
			\hat{P}\left(\frac{k}{n/l}, \yif{k}\right) = F + k,
		\]
		and
		\[
			D = n/l - \frac{1}{2}(kl-n) - \yif{k}.
		\]
	\end{lem}
	\begin{proof}
        All values on the same staircase $\hat{J}_i$ belong to the same congruence class
        modulo $k$. More specifically, since $k > 0$ by assumption, the values on the
        $\size{\hat{J}_i}$ 
        steps are consecutive mod $k$-values starting with the smallest integer in the congruence class
        that is at least $\hat{P}\left(\frac{i}{n/l},0\right)$.
        For each $i$, let $c(i)$ be the smallest positive integer such that the values of
        $\hat{J}_{i+c(i)}$ are in the same congruence class as $\hat{J}_i$.
        From a certain value $i_0$, we have
        $\hat{P}\left(\frac{i}{n/l}, 0\right) < \hat{P}\left(\frac{i+c(i)}{n/l}, 0\right)$
        for all $i > i_0$.
        This means that the value of $\hat{P}$ on the first step of $\hat{J}_{i+c(i)}$
        must be greater than the value of $\hat{P}$ on the last step of $\hat{J}_i$,
        since it would otherwise have to equal one of the values on $\hat{J}_i$.
    
        For large enough $i$, the value of $\hat{P}$ on the last step of $\hat{J}_i$
        must be exactly $k$ less than the value on the first step of $\hat{J}_{i+c(i)}$.
        All non-negative values in the congruence class must be represented by $\hat{P}$ exactly once.
        Since the values are strictly increasing for $i > i_0$, if there is a gap between the last
        step of $\hat{J}_i$ and the first step of $\hat{J}_{i+c(i)}$, then the ``missing'' values
        must be found on staircases to the left of $\hat{J}_{i_0}$.
        There are finitely many such staircases, so from a certain point, there can be no gaps.
        
        The points on the staircases $\hat{J}_{(n/l)i}$ have integral $x$-coordinate.
        The first step is $(i, 0)$, the last step is $ni$, and there are $li + 1$ steps,
        cf.~Remark \ref{rem:trfSectorStructure}.
        We further note that $\yif{(n/l)i + j} = \yif{j}$ for all $i,j$.
        For all large enough $i$, we can thus write
		\begin{equation}\label{eq:Pini1}
			\hat{P}(i, ni) = \hat{P}\left(i + \frac{c(i)}{n/l}, \yif{c(i)} \right) - k.
		\end{equation}
		Since $\size{\hat{J}_{(n/l)i}} = li + 1$, we can also write
		\begin{equation}\label{eq:Pini2}
		    \hat{P}(i, ni) = \hat{P}(i,0) + k l i.
		\end{equation}
		Using formula \eqref{eq:initNecForm}, we get
		\begin{align*}
		    \hat{P}\left(\!i + \frac{c(i)}{n/l}, \yif{c(i)} \!\right)
		        &= \frac{n}{2}\left(\!i + \frac{c(i)}{n/l}\!\right)\!
		            \left(\!i + \frac{c(i)}{n/l} - 1\!\right)
		            + D\left(\!i + \frac{c(i)}{n/l}\!\right) + \frac{k}{n/l} \yif{c(i)} + F\\
		        &= \frac{n}{2} i (i-1) + Di + F\\
		            &\hspace{1cm}+ \frac{n}{2} \frac{c(i)}{n/l}\left(\frac{c(i)}{n/l}-1\right)
		                +D \frac{c(i)}{n/l} + \frac{k}{n/l} \yif{c(i)} + F\\
		            &\hspace{1cm}-F + c(i)l i\\
		        &= \hat{P}(i,0) + \hat{P}\left(\frac{c(i)}{n/l}, \yif{c(i)}\right) - F +c(i)li.
		\end{align*}
		Combined with \eqref{eq:Pini1} and \eqref{eq:Pini2}, it follows that
		\[
		    \hat{P}\left(\frac{c(i)}{n/l}, \yif{c(i)}\right) = k + F + (k-c(i))li.
		\]
		This holds for all $i$ large enough. For only finitely many $i$ can we have $c(i) < k$,
		as the formula would otherwise produce more than one value for
		$\hat{P}\left(\frac{c(i)}{n/l}, \yif{c(i)}\right)$.
		Also, we cannot have $c(i) > k$ for infinitely many $i$ as we would then eventually
		have a negative value for $\hat{P}\left(\frac{c(i)}{n/l}, \yif{c(i)}\right)$.
		We can therefore find an $i$ large enough that the formula must hold and for which
		$c(i) = k$. Using this in the formula, we obtain
		\[
		    \hat{P}\left(\frac{k}{n/l}, \yif{k}\right) = k + F.
		\]
\begin{figure}
	\centering
    \includegraphics{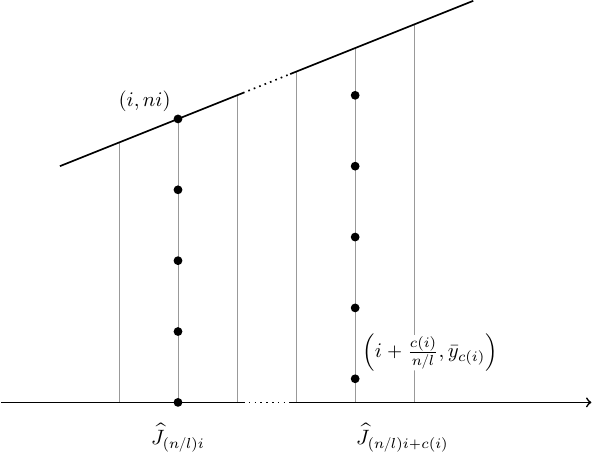}
	\caption{Next staircase with the same congruence class mod $k$
	after a staircase with integral $x$-coordinate.
	$\hat{P}(i, ni) \equiv \hat{P}\left(i + \frac{c(i)}{n/l}, \yif{c(i)}\right) \pmod k$.}
	\label{fig:nextModkCongCl}
\end{figure}

		Expanding $\hat{P}\left(\frac{k}{n/l}, \yif{k} \right)$ leads to
		\[
			D = n/l - \frac{1}{2}(kl - n) - \yif{k}.
		\]
	\end{proof}

    \begin{example}\label{ex:throughEx2}
        The polynomial in Example \ref{ex:throughEx1},
        in the shape of \eqref{eq:initNecForm}, 
        \[
            \hat{P}(x,y) = 6x \left(x-1\right) - 2 x + \frac{3}{2} y + 2.
        \]
        is a QPP on $\hat{I}\left(\frac{12}{7}\right)$ with $k=3$.
        We have $l = 6$, $n/l = 2$, and $\yif{k} = 1$ (see Figure \ref{fig:12o7kp3Ex1}).
        We verify that
        \[
            \hat{P}\left(\frac{3}{2}, 1\right) = 5 = k + F,
        \]
        and
        \[
            D = -2 = 2 - \frac{1}{2}(3\cdot 6 - 12) - 1.
        \]
    \end{example}
	
	\begin{lem}\label{lem:modkStairs}
	    Let $\hat{P}(x,y)$ be a QPP on $\hat{I}\left(\frac{n}{m}\right)$ with $k > 0$.
		Then for all $i$,
		the values of $\hat{P}$ on $\hat{J}_{i+k}$ are in the same congruence class modulo $k$
		as the values on $\hat{J}_i$.
    \end{lem}
	\begin{proof}
		For any $i$,
		pick $\left(\frac{i}{n/l}, \yif{i}\right)$ and $\left(\frac{i+k}{n/l}, \yif{i+k}\right)$
		as representatives of points from $\hat{J}_i$ and $\hat{J}_{i+k}$, respectively.
        By formula \eqref{eq:initNecForm}, we have
        \begin{align*}
            \hat{P}\left(\!\frac{i+k}{n/l}, \yif{i\!+\!k} \!\right) &=
                    \frac{n}{2}\!\left(\!\frac{i}{n/l} \!+\! \frac{k}{n/l}\!\right)\!
                        \left(\!\frac{i}{n/l} \!+\! \frac{k}{n/l} - 1\!\right)
                    + D\!\left(\!\frac{i}{n/l}\!+\!\frac{k}{n/l}\!\right)
                    + \frac{k}{n/l}\yif{i\!+\!k} + F\\[0.2cm]
                &= \frac{n}{2} \frac{i}{n/l}\left(\frac{i}{n/l} -1\right)
                            +D\frac{i}{n/l} + \frac{k}{n/l} \yif{i} + F\\
                    &\hspace{1cm}+  \frac{n}{2} \frac{k}{n/l}\left(\frac{k}{n/l} -1\right)
                            +D\frac{k}{n/l} + \frac{k}{n/l} \yif{k} + F\\
                    &\hspace{1cm}+ ik l^2/n + \frac{k}{n/l} \yif{i+k}
                            - \frac{k}{n/l} \yif{i} - \frac{k}{n/l}\yif{k} - F\\[0.2cm]
                &= \hat{P}\left(\!\frac{i}{n/l}, \yif{i} \!\right)
                            + \hat{P}\left(\!\frac{k}{n/l}, \yif{k} \!\right) - F
                            + k\!\left(\!i l^2/n + \frac{\yif{i+k} \!-\! (\yif{i} \!+\! \yif{k})}{n/l}\!\right).
        \end{align*}
        By Lemma \ref{lem:nextCong}, $\hat{P}\left(\frac{k}{n/l}, \yif{k} \right) = k + F$, so
        \begin{equation}\label{eq:equivModk}
            \hat{P}\left(\frac{i+k}{n/l}, \yif{i+k} \right) =
                \hat{P}\left(\frac{i}{n/l}, \yif{i} \right) + k
                        + k\left(i l^2/n + \frac{\yif{i+k} - (\yif{i} + \yif{k})}{n/l}\right).
        \end{equation}
		In the integral sector case, all stairs begin at $0$.
		Else, cf.~Remark \ref{rem:trfSectorStructure}, we have
		\[
			\yif{i+k} \equiv -\frac{1}{(m-1)/l}(i+k)  \equiv \yif{i} + \yif{k} \pmod {n/l},
		\]
		so $\frac{\yif{i+k} - (\yif{i} + \yif{k})}{n/l}$ is an integer
		and \eqref{eq:equivModk} shows that
		\[
		    \hat{P}\left(\frac{i+k}{n/l}, \yif{i+k} \right) \equiv
		        \hat{P}\left(\frac{i}{n/l}, \yif{i} \right) \pmod k
		\]
		which is what we wanted to show.
	\end{proof}

    \begin{example}\label{ex:throughEx3}
        Note, for the polynomial
        $\hat{P}(x,y) = 6x \left(x-\frac{3}{2}\right) + x + \frac{3}{2} y + 2$
        from Examples \ref{ex:throughEx1} and \ref{ex:throughEx2},
        which is a QPP with $k = 3$,
        that the values on any two staircases that are three apart
        are congruent modulo $3$. See Figure \ref{fig:12o7kp3Ex1}.
    \end{example}
	
	\begin{lem}\label{lem:ycoordk}
	    Let $\hat{P}(x,y)$ be a QPP on $\hat{I}\left(\frac{n}{m}\right)$ with $k > 0$.
		Then
		\[
			\yif{k} = n/l - 1.
		\]
	\end{lem}
	\noindent Note that this is trivially true for integral sectors.
	\begin{proof}
		For large enough $i$, the value of $\hat{P}$ on the first step of
		$\hat{J}_i$ is $k$ more than the value on the last step of $\hat{J}_{i-k}$.
        This follows from Lemma \ref{lem:modkStairs}, because $\hat{P}$ is eventually increasing
        in the first coordinate, cf.~the argument given in the proof of Lemma \ref{lem:nextCong}.
    	The line $x = \frac{i-k}{n/l}$ intersects the line $y = nx$ at the $y$-coordinate $l(i-k)$
		(see Figure \ref{fig:firstStepOnStairk}),
		so, since $\hat{P}$ is strictly increasing in the $y$-coordinate, we have
		\begin{equation}\label{eq:ineq}
			\hat{P}\left(\frac{i-k}{n/l}, l(i-k)\right) + k \geq \hat{P}\left(\frac{i}{n/l}, \yif{i}\right).
		\end{equation}
        \begin{figure}[hbtp]
        	\centering
            \includegraphics{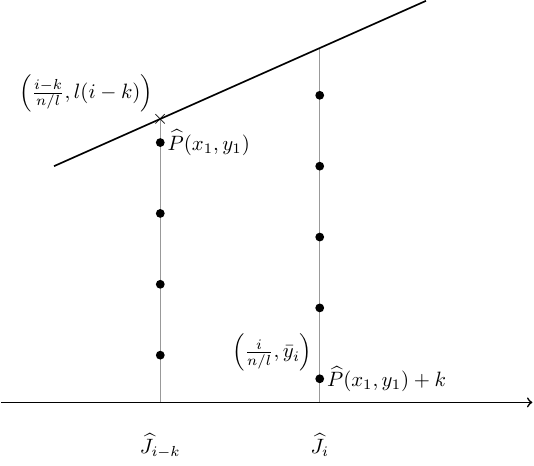}
        	\caption{The value of $\hat{P}(x,y)$ is $k$ more on the first step of $\hat{J}_i$
        	than on the last step of $\hat{J}_{i-k}$.
        	The lines $x = \frac{i-k}{n/l}$ and $y = nx$ intersect at
        	$\left(\frac{i-k}{n/l}, l(i-k)\right)$.}
        	\label{fig:firstStepOnStairk}
        \end{figure}
        By Formula \eqref{eq:initNecForm}, we have
        \begin{align*}
            \hat{P}\left(\frac{i-k}{n/l}, \yif{i-k}\right)
                &= \frac{n}{2} \frac{i-k}{n/l} \left(\frac{i-k}{n/l}-1\right)
                    + D \frac{i-k}{n/l} + \frac{k}{n/l}l(i-k) + F\\
                &= \frac{n}{2}\frac{i}{n/l}\left(\frac{i}{n/l}-1\right)
                            + D\frac{i}{n/l} + \frac{k}{n/l}\yif{i} + F\\
                    &\hspace{1cm} -\frac{1}{2}k^2l^2/n + \frac{1}{2}kl
                            - D\frac{k}{n/l} - \frac{k}{n/l}\yif{i}\\
                &= \hat{P}\left(\frac{i}{n/l}, \yif{i}\right)
                    -\frac{1}{2}k^2l^2/n + \frac{1}{2}kl
                            - D\frac{k}{n/l} - \frac{k}{n/l}\yif{i}.
        \end{align*}
        By Lemma \ref{lem:nextCong},
        $D = n/l - \frac{1}{2}(kl-n) - \yif{k}$, so the above equation simplifies to
		\[
			\hat{P}\left(\frac{i-k}{n/l}, l(i-k)\right) + k - \hat{P}\left(\frac{i}{n/l}, \yif{i}\right)
			= \frac{k}{n/l}(\yif{k} - \yif{i}).
		\]
		The left-hand side is non-negative by \eqref{eq:ineq}.
		Since this holds true for all $i$, the value of $\yif{k}$ must be the largest possible,
		which is $n/l - 1$.
	\end{proof}

    \begin{example}\label{ex:throughEx4}
        For our recurring QPP with $k=3$ on the sector 
        $I\left(\frac{12}{7}\right)$ (Examples \ref{ex:throughEx1},
        \ref{ex:throughEx2} and \ref{ex:throughEx3}),
        we have $\yif{k} = n/l - 1 = 1$. See Figure \ref{fig:12o7kp3Ex1}.
    \end{example}
	
	\begin{prop}\label{prop:necCond}
		If $\hat{P}(x,y)$ is a QPP on $\hat{I}\left(\frac{n}{m}\right)$ with $k > 0$ and $l = \gcd(n,m-1)$,
		then we must have $n/l \mid l$ and
		\begin{equation}\label{eq:necForm}
			\hat{P}(x,y) = \frac{n}{2} x \left( x - \frac{k}{n/l}\right) + x + \frac{k}{n/l} y + F,
		\end{equation}
		where $k \equiv (m-1)/l \pmod {n/l}$.
	\end{prop}
	\begin{proof}
        The condition $n/l \mid l$ follows from the necessary condition $n \mid (m-1)^2$
        of Theorem \ref{thm:stanNecForm}. 
        Since $n \mid (m-1)^2$ implies $n/l \mid (m-1)^2/l$,
        and because $n/l$ and $(m-1)/l$ are coprime, we must have $n/l \mid l$.
        
		The necessary form of $\hat{P}(x,y)$ follows directly from the initial form in Equation \eqref{eq:initNecForm}
		and Lemmas \ref{lem:nextCong} and \ref{lem:ycoordk}.
		The condition $k \equiv (m-1)/l \pmod {n/l}$ follows from Lemma~\ref{lem:ycoordk} and
		the fact that $\yif{i} (m-1)/l  \equiv - i \pmod {n/l}$.
		Again, this is trivial for integral sectors.
	\end{proof}

\section{Classification}

In this section we prove our main result, Theorem \ref{thm:main}.
First, we give two more lemmas. Lemma \ref{lem:sufCond} gives 
a sufficient condition for $\hat{P}$ to be a packing polynomial. 
Lemma \ref{lem:Fform} gives the constant term $F$.

	\begin{lem}\label{lem:sufCond}
	    Let $m,n \in \N$, $l = \gcd(m-1, n)$ and
	    \[
	        \hat{P}(x,y) = \frac{n}{2} x \left( x - \frac{k}{n/l}\right) + x + \frac{k}{n/l} y + F,
	    \]
	    where $k \equiv (m-1)/l \pmod {n/l}$.
		$\hat{P}(x,y)$ is a QPP with $k > 0$ on $\hat{I}\left(\frac{n}{m}\right)$ if and only if
		\[
			\left\{ \hat{P}\left(\frac{i}{n/l}, \yif{i}\right) \forwhich i = 0, 1, \ldots, k-1 \right\}
			= \{0,1,\ldots,k-1\}.
		\]
	\end{lem}
	\begin{proof}
		By Lemma~\ref{lem:modkStairs}, for all $i$ we have that values of
        $\hat{P}$ on $\hat{J}_{i+k}$ are in the same congruence class modulo $k$
        as the values of $\hat{P}$ on $\hat{J}_{i}$.
		Therefore, if $\hat{P}(x,y)$ takes the values $0, 1, \ldots, k-1$ (in any order)
		on the first steps of the first $k$ staircases,
		then $\hat{P}(x,y)$ is a QPP on $\hat{I}\left(\frac{n}{m}\right)$
		if
		\begin{equation}\label{eq:difNextStair}
			\hat{P}\left(\frac{i+k}{n/l}, \yif{i+k}\right) -\hat{P}\left(\frac{i}{n/l}, \yif{i}\right)
			= k \left(\size{\hat{J_i}}\right).
		\end{equation}
        Conversely, if \eqref{eq:difNextStair} holds and
        $\hat{P}(x,y)$ fails to take the values $0, 1, \ldots, k-1$
		on the first steps of the first $k$ staircases,
        then, since $\hat{P}$ is increasing in both coordinates,
        any missing value could never be attained again on $\hat{I}\left(\frac{n}{m}\right)$.
        Therefore $\hat{P}$ would not be surjective onto $\N_0$ and thus not a QPP.

		Applying formula \eqref{eq:necForm}, a direct computation shows
		\[
			\hat{P}\left(\frac{i+k}{n/l}, \yif{i+k}\right) 
                -\hat{P}\left(\frac{i}{n/l}, \yif{i}\right)
			= k \left( (l^2/n)i + \frac{1 + \yif{i+k} - \yif{i}}{n/l} \right).
		\]
		Since $\yif{i+k} - \yif{i} \equiv -1 \pmod {n/l}$ we have
		$\yif{i+k} - \yif{i} = (n/l)\llbracket n/l \mid i \rrbracket - 1$
		and the result follows.
	\end{proof}

    \begin{example}\label{ex:throughEx5}
        For the recurring example,
        $\hat{P}(x,y) = 6x \left(x-\frac{3}{2}\right) + x + \frac{3}{2} y + 2$,
        a QPP on $\hat{I}\left(\frac{12}{7}\right)$ with $k = 3$,
        the first steps on the first three staircases are
        $\yif{0} = 0$, $\yif{1} = 1$, $\yif{2} = 0$, and we can verify that
        \begin{align*}
            \hat{P}\left(0,0\right) &= 2, \quad
            \hat{P}\left(\frac{1}{2},1\right) = 1, \quad
            \hat{P}\left(1,0\right) = 0.
        \end{align*}
        We can now with confidence say that it indeed is a QPP on
        $\hat{I}\left(\frac{12}{7}\right)$.
    \end{example}
	
	\begin{lem}\label{lem:Fform}
		If $\hat{P}(x,y)$ is a QPP with $k > 0$ on $\hat{I}\left(\frac{n}{m}\right)$, then
		\begin{equation}\label{eq:Fform}
			F = \frac{(l^2/n) (k-1)(k + 1)}{12}.
		\end{equation}
	\end{lem}
	\begin{proof}
		By Lemma~\ref{lem:sufCond}, we have
		\[
			\sum_{i=0}^{k-1} \hat{P}\left(\frac{i}{n/l}, \yif{i} \right) = \sum_{i=0}^{k-1} i,
		\]
		from which we get
		\[
			\frac{l^2/n}{2} \sum_{i=0}^{k-1} i^2
			+ \left(\frac{l}{n}-\frac{k}{2}l^2/n - 1\right) \sum_{i=0}^{k-1} i
			+ \frac{k}{n/l} \sum_{i=0}^{k-1} \yif{i} + k F = 0.
		\]
		Observing that $\yif{i} + \yif{k-i} \equiv -1 \pmod {n/l}$ which implies that
		$\yif{i} + \yif{k-i} = n/l - 1$, we obtain,
		remembering that $\yif{0} = 0$,
		by rearranging the terms the formula
		\[
			\sum_{i=0}^{k-1} \yif{i} =
			\left\{\begin{array}{ll}\displaystyle
						\sum_{i=1}^{(k-1)/2} (\yif{i} + \yif{k-1}), & \il{$k$ odd}\\[0.5cm]\displaystyle
						\sum_{i=1}^{k/2 - 1} (\yif{i} + \yif{k-1}) + \yif{k/2}, & \il{$k$ even}
			       \end{array}
			\right\} = \frac{(k-1)(n/l-1)}{2}.
		\]
		The formula for $F$ now follows by applying the well-known
		\[
			\sum_{i=0}^{k-1} i^2 = \frac{(k-1)k(2k - 1)}{6} \il{ and }
			\sum_{i=0}^{k-1} i = \frac{(k-1)k}{2}.
		\]
	\end{proof}

    \begin{example}\label{ex:throughEx6}
        As we have seen (Examples \ref{ex:throughEx1}, \ref{ex:throughEx2},
        \ref{ex:throughEx3}, \ref{ex:throughEx4} and \ref{ex:throughEx5}),
        The polynomial 
        $\hat{P}(x,y) = 6x \left(x-\frac{3}{2}\right) + x + \frac{3}{2} y + 2$
        is a QPP on $\hat{I}\left(\frac{12}{7}\right)$ with $k = 3$.
        We have $n= 12$, $l = 6$, so $l^2/n = 3$ and we can verify that
        \[
            F = \frac{3 (3-1)(3+1)}{12} = 2.
        \]
    \end{example}

	\begin{thm}
	\label{thm:main}
		Let $m,n\in \N$ with $\gcd(m,n) = 1$ and put $l = \gcd(m-1, n)$.
		If $n \mid l^2$ and
		\[\renewcommand{\arraystretch}{1.5}
    	   	\begin{array}{llll}
    		    (1) & (m-1)/l &\equiv \pm 1 \pmod {n/l}, &\il{set } k = \pm 1, \il{ or }\\
    		    (2) & (m-1)/l &\equiv \pm 2 \pmod {n/l} \il{ and } l^2/n = 4, &\il{set } k = \pm 2, \il{ or }\\
    		    (3) & (m-1)/l &\equiv \pm 3 \pmod {n/l} \il{ and } l^2/n = 3, &\il{set } k = \pm 3,
    		\end{array}
		\]
		then
		\[
			P(x,y) = \frac{n}{2} \left(x - \frac{m-1}{n}y\right)\left(x - \frac{m-1}{n}y - \frac{kl}{n}\right)
				+ x + \frac{kl - (m-1)}{n} y + |k| - 1
		\]
		is a QPP on $I\left(\frac{n}{m}\right)$.
		These are the only QPPs on sectors $S(\alpha) \subseteq \R^2$.
	\end{thm}
	\begin{proof}
		Assume that $n/l \mid l$ and let $\hat{P}(x,y)$ have the necessary form in Equation \eqref{eq:necForm}.
		Assume first that $k > 0$. We go over the cases $k = 1, 2, 3,$ and $k > 4$ individually.
		By Proposition \ref{prop:necCond}, it is in each case a requirement that
		$(m-1)/l \equiv k \pmod {n/l}$ (again trivial in the integral case, as $n/l = 1$).
        We determine whether
        $\hat{P}$ takes the values $0, 1, \ldots, k-1$ (in any order) on the first steps
        of the first $k$ staircases.
        If that is the case, then Lemma \ref{lem:sufCond} guarantees a packing polynomial.
        
		\begin{case}{$k=1$}
			By Lemma \ref{lem:Fform}, we must have $F = 0$. Since
			\[
				\hat{P}(0,0) = 0,
			\]
			the requirements of Lemma \ref{lem:sufCond} are met.
		\end{case}
		
		\begin{case}{$k=2$}
			By Lemma \ref{lem:sufCond}, we must have $F = 0$ or $1$. Since $l^2/n > 0$,
			we must have $l^2/n = 4$ and $F = 1$, by Lemma~\ref{lem:Fform}.
			Since $(m-1)/l\equiv 2\pmod {n/l}$, $n/l$ is odd and we have
			in the non-integral sector case
			$\yif{1} \equiv - \frac{1}{(m-1)/l} \equiv -\frac{1}{2} \pmod {n/l}$,
			so $\yif{1} = \frac{n/l - 1}{2}$. We find that
			\[
				\hat{P}(0,0) = 0
			\]
			and
			\[
				\hat{P}\left(\frac{1}{n/l}, \frac{n/l - 1}{2}\right) = 1.
			\]
			By Lemma~\ref{lem:sufCond}, $\hat{P}(x,y)$ is a QPP on $\hat{I}\left(\frac{n}{m}\right)$.
		\end{case}
		
		\begin{case}{$k=3$}
			We must have $F = 0, 1$ or $2$. By Lemma~\ref{lem:Fform}, we must have
			$F = \frac{2 l^2/n}{3} = 1$ or $2$. Since this must be an integer, we conclude that
			$F = 2$ and $l^2/n = 3$.
			Since $(m-1)/l\equiv 3\pmod {n/l}$, we have $n/l \not\equiv 0 \pmod 3$ and
			$3\yif{1} \equiv -1 \pmod {n/l}$ and $3\yif{2} \equiv -2 \pmod {n/l}$.
			This means that
			\[
				\yif{1} = \left\{\begin{array}{ll}
									\frac{n/l - 1}{3} & \il{if } n/l \equiv 1 \pmod 3\\
									\frac{2 n/l - 1}{3} & \il{if } n/l \equiv 2 \pmod 3,
				                 \end{array}\right.
			\]
			and
			\[
				\yif{2} = \left\{\begin{array}{ll}
									\frac{2(n/l - 1)}{3} & \il{if } n/l \equiv 1 \pmod 3\\
									\frac{n/l - 2}{3} & \il{if } n/l \equiv 2 \pmod 3.
				                 \end{array}\right.
			\]
			In either case, we have
			\[
				\hat{P}(0,0) = 2 \il{ and }
				\left\{
					\hat{P}\left (\frac{1}{n/l}, \yif{1} \right), \hat{P}\left(\frac{2}{n/l}, \yif{2}\right)
				\right\}
				= \{0,1\}.
			\]
			Again, Lemma~\ref{lem:sufCond} verifies that $\hat{P}(x,y)$ is a QPP on $\hat{I}\left(\frac{n}{m}\right)$.
		\end{case}
		
		\begin{case}{$k > 3$}
			By Lemma~\ref{lem:Fform},
			we must have $l^2/n (k+1) \leq 12$, since $0 \leq F \leq k-1$ by Lemma~\ref{lem:sufCond}.
			This means that $4 \leq k \leq 11$. For $k = 4$ and $5$, this implies that $l^2/n = 1$ or $2$,
			and for $k > 5$ that $l^2/n = 1$.
			Since $F$ is an integer, checking each possible value for $k$ and $l^2/n$ in the formula
			\eqref{eq:Fform} leaves as the only options $k = 5, 7$ or $11$.
			In particular, $k$ must be odd.
			Since $k \equiv (m-1)/l \pmod {n/l}$, this means that both $2$, $l^2/n$ and $n/l$ have
			inverses modulo $k$ and thus likewise for $l = (l^2/n)(n/l)$ and $n$.
			Therefore
			\begin{align*}
				\hat{P}\left (\frac{i}{n/l}, \yif{i} \right)
					&= \frac{n}{2}\frac{i}{n/l}\left(\frac{i}{n/l} - \frac{k}{n/l}\right)
						+ \frac{i}{n/l} + \frac{k}{n/l} \yif{i} + F\\
					&\equiv \frac{l^2/n}{2} i^2 + \frac{i}{n/l} + F \pmod k\\
					&\equiv \frac{l^2/n}{2} \left(i + \frac{1}{l}\right)^2 - \frac{1}{2n} + F \pmod k.
			\end{align*}
			For any $j \equiv -\frac{1}{l} \pmod {k}$, we have
			\[
				\hat{P}\left (\frac{j+1}{n/l}, \yif{j + 1} \right)
				\equiv
				\hat{P}\left (\frac{j-1}{n/l}, \yif{j - 1} \right) \pmod k.
			\]
			This is impossible by Lemma~\ref{lem:modkStairs}.
		\end{case}
		
		We have now determined all possibilities for QPPs with $k > 0$.
		We will deal with the negative cases by showing that there is a one-to-one
		correspondence between QPPs with $k = k'$ and QPPs with $k = -k'$.
		
		Assume that $\hat{P}(x,y)$ is a QPP on $\hat{I}\left(\frac{n}{m}\right)$ with $k < 0$.
        By Theorem \ref{thm:stanNecForm} and Equation \eqref{eq:initNecForm}, we must have $n/l\mid l$ and $\hat{P}(x,y)$ must have the form
		\[
			\hat{P}(x,y) = \frac{n}{2} x (x-1) + Dx + \frac{k}{n/l} y + F
		\]
		with $D, F \in \Z$. The involutory transformation $L:(x,y) \mapsto (x, nx - y)$
		maps $S(n)$ onto $S(n)$ while vertically flipping the coordinates.
		See Figure \ref{fig:ascVsDesc}.
\begin{figure}
	\centering
    \includegraphics{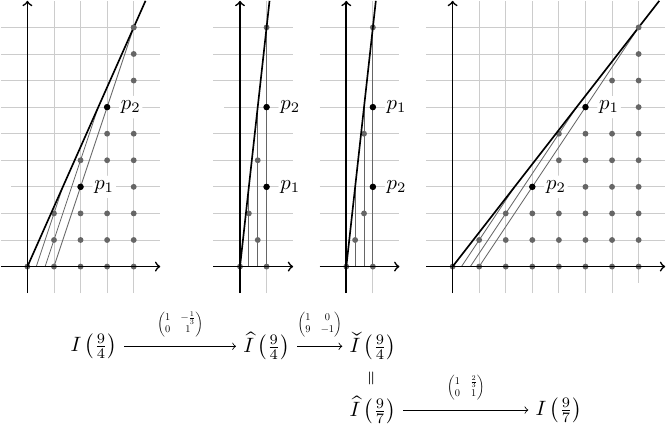}
	\caption{A QPP with $k$ negative corresponds to a QPP with $k$ positive.
	    The sectors are not necessarily the same.}
	\label{fig:ascVsDesc}
\end{figure}
		Specifically,
		the effect on the points is
		\begin{align*}
			\check{I}\left(\frac{n}{m}\right) &:= L\left(\hat{I}\left(\frac{n}{m}\right)\right)\\
				&= \left\{
					(x,y) \in S(n) \forwhich x =\frac{i}{n/l}, ((m-1)/l) y \equiv i \pmod {n/l},
					i \in \N_0
				\right\}.
		\end{align*}
		Then
		\[
			\check{P}(x,y) = \hat{P}(x,nx-y)
				= \frac{n}{2} x(x-1) + (D + kl)x - \frac{k}{n/l} y + F
		\]
		is a QPP on $\check{I}\left(\frac{n}{m}\right)$ and the difference in values of consecutive
		steps is $-k$.
		From the above analysis, we have
		\[
			\yif{-k} = n/l - 1 \equiv \frac{-k}{(m-1)/l} \pmod {n/l}, \il{ so }
			(m-1)/l \equiv k \pmod {n/l}
		\]
		and
		\[
			\check{P}(x,y) = \frac{n}{2} x \left(x-\frac{-k}{n/l}\right) + x + \frac{-k}{n/l} + F,
		\]
		where either
		\begin{align*}
		                -k &= 1 \il{ and } F = 0,\\
		    \il{or }    -k &= 2,~  l^2/n = 4, \il{ and } F = 1,\\
		    \il{or }    -k &= 3,~  l^2/n = 3, \il{ and } F = 2.
		\end{align*}
		This means that
		\[
			\hat{P}(x,y) = \check{P}(x, nx - y) =
			\frac{n}{2} x \left(x-\frac{k}{n/l}\right) + x + \frac{k}{n/l} y + (-k) - 1.
		\]
		So, regardless of the sign of $k$, if $k = \pm 1, \pm 2, \pm 3 \equiv (m-1)/l \pmod {n/l}$,
		then
		\[
			\hat{P}(x,y) = \frac{n}{2} x \left(x-\frac{k}{n/l}\right) + x + \frac{k}{n/l} y + |k| - 1
		\]
		is a QPP on $\hat{I}\left(\frac{n}{m}\right)$
		and these are the only options.
		Transforming back to the original integral sector, we find that
		\[
			P(x,y) = \frac{n}{2} \left(x - \frac{m-1}{n}y\right)\left(x - \frac{m-1}{n}y - \frac{kl}{n}\right)
				+ x + \frac{kl - (m-1)}{n} y + |k| - 1
		\]
		is a QPP on $I\left(\frac{n}{m}\right)$.
	\end{proof}
	
    \begin{example} 
        \label{ex:8over5}
        Let $\alpha = \frac{8}{5}$. Then we have
        $n = 8$, $m = 5$, so $n \mid (m-1)^2$. Then
        $l = \gcd(m-1, n) = 4$, $n/l = 2$ and
        $
            (m-1)/l = 1 \equiv \pm 1, \pm 3 \pmod 2,
        $
        but $l^2/n = 2 \neq 3$.
        So only $k = \pm 1$ yields QPPs. That is, the polynomials
        \begin{align*}
            P(x,y) &= 4\left(x - \frac{1}{2} y \right)
                        \left(x - \frac{1}{2} y - \frac{1}{2}\right) + x,\\
            P(x,y) &= 4\left(x - \frac{1}{2} y \right)
                        \left(x - \frac{1}{2} y + \frac{1}{2}\right) + x - y
        \end{align*}
        are the only QPPs on $I\left(\frac{8}{5}\right)$, see Figure \ref{fig:8over5}.
        \begin{figure}[htbp]
            \centering
            \includegraphics{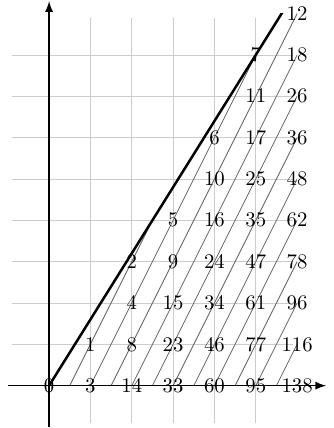}\hspace{0.5cm}
            \includegraphics{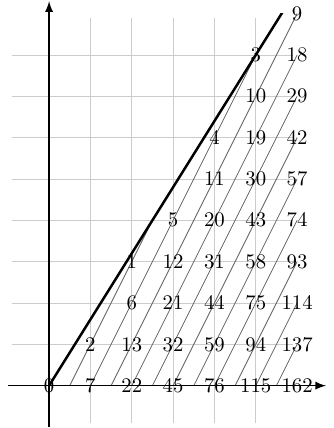}
            \caption{The quadratic packing polynomials on $S(8/5)$ given in Example \ref{ex:8over5}.}
            \label{fig:8over5}
        \end{figure}
    \end{example}

    \begin{example} 
        \label{ex:12over7}
        Let $\alpha = \frac{12}{7}$.
        Then we have $n = 12$, $m = 7$, so $n \mid (m-1)^2$.
        Then $l = \gcd(m-1, n) = 6$, $n/l = 2$ and
        $
            (m-1)/l = 1 \equiv \pm 1, \pm 3 \pmod 2,
        $
        and $l^2/n = 3$.
        So setting $k = \pm 1, \pm 3$ yields QPPs. That is, the polynomials
        \begin{align*}
            P(x,y) &= 6\left(x - \frac{1}{2} y \right)
                        \left(x - \frac{1}{2} y - \frac{1}{2}\right) + x,\\
            P(x,y) &= 6\left(x - \frac{1}{2} y \right)
                        \left(x - \frac{1}{2} y + \frac{1}{2}\right) + x - y\\
            P(x,y) &= 6\left(x - \frac{1}{2} y \right)
                        \left(x - \frac{1}{2} y - \frac{3}{2}\right) + x + y + 2\\
            P(x,y) &= 6\left(x - \frac{1}{2} y \right)
                        \left(x - \frac{1}{2} y + \frac{3}{2}\right) + x - 2 y + 2
        \end{align*}
        are the only QPPs on $I\left(\frac{12}{7}\right)$, see Figure \ref{fig:12over7}. In fact, $S(3)$, $S(4)$ and $S\left(\frac{12}{7}\right)$ are the only sectors (up to equivalence of transformation) with four quadratic packing polynomials. 
        \begin{figure}[htbp]
            \centering
            \includegraphics{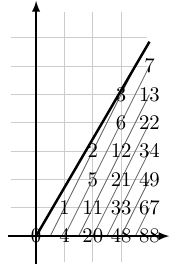}\hspace{-2pt}
            \includegraphics{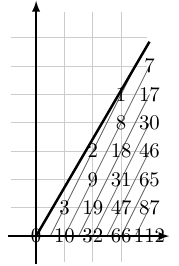}\hspace{-2pt}
            \includegraphics{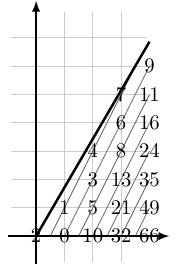}\hspace{-2pt}
            \includegraphics{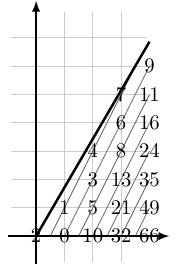}\\
            $k = 1$\hspace{2.1cm} $k= -1$\hspace{1.8cm} $k = 3$\hspace{2.1cm} $k= -3$
            \caption{The quadratic packing polynomials on $S(12/7)$ given in Example \ref{ex:12over7}.}
            \label{fig:12over7}
        \end{figure}
    \end{example}
    
    \begin{example}
        \label{ex:n}
        Consider the integral sectors.
        That is, $n \in \N$, $m = 1$.
        Then $n \mid 0$, 
        $l = \gcd(0, n) = n$, $n/l = 1$ and
        $
            (m-1)/l = 0 \equiv \pm 1, \pm 2, \pm 3 \pmod 1.
        $
        So setting $k = \pm 1$ yields QPPs for all $n$.
        As $l^2/n = n$, for $n = 3$ and $n = 4$ we obtain the
        two additional QPPs discovered by Stanton \cite{stan14}.
        
        So for all $n$, the polynomials \eqref{eq:nathintpols} discovered by Nathanson \cite{nath14} are QPPs.
        Additionally, for $n = 3$ the polynomials
        \begin{align*}
            P(x,y) &= \frac{3}{2} x (x \pm 3) + x \mp 3 y + 2
        \end{align*}
        and for $n = 4$ the polynomials
        \begin{align*}
            P(x,y) &= 2 x (x \pm 2) + x \mp 2 y + 1
        \end{align*}
        are QPPs, see Figure \ref{fig:n}.
        
        \begin{figure}[htbp]
            \centering
            \includegraphics{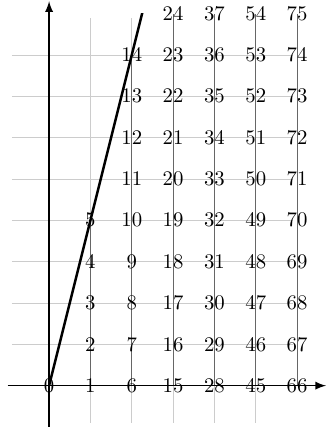}\hspace{0.5cm}
            \includegraphics{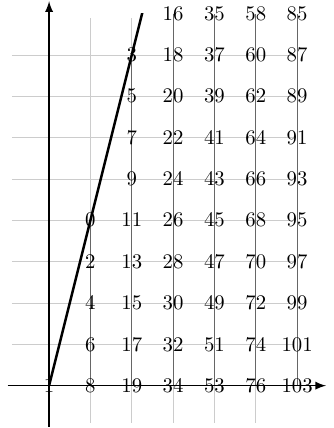}\\
            $k=1$\hspace{5cm} $k=-2$
            \caption{Two quadratic packing polynomials on $S(4)$ given in Example \ref{ex:n}.}
            \label{fig:n}
        \end{figure}
    \end{example}
    
    \begin{rem}
        \label{ex:cantor}
        While $S(\infty)$ is not explicitly included by Theorem \ref{thm:main} as it is stated,
        it does follow as a corollary since the sector is equivalent to $S(1)$, 
        cf.~ Remark \ref{rem:trfSectorInfty}.
        For $S(\infty)$, we set $n=1$, $m=0$, and note that indeed $1 \mid (-1)^2$.
        We have $l = \gcd(-1, 1) = 1$, $n/l = 1$ and
        $
            (m-1)/l = -2 \equiv \pm 1, \pm 2, \pm 3 \pmod 1.
        $
        So setting $k = \pm 1$ yields QPPs, whereas $l^2/n \neq 3, 4$.
        We recover the original Cantor polynomials in Equation \eqref{eq:cantor}, see Figure \ref{fig:cantor}.
        \begin{figure}[htbp]
            \centering
            \includegraphics{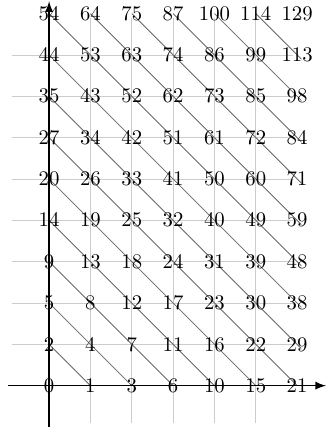}\hspace{0.5cm}
            \includegraphics{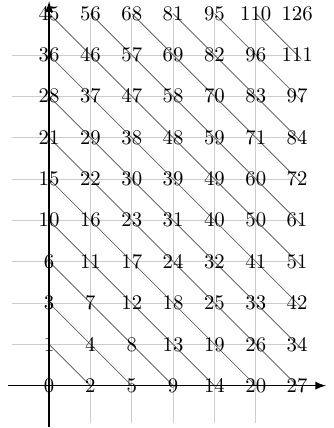}\\
            $k=1$\hspace{5cm} $k=-1$
            \caption{The Cantor polynomials on $S(\infty)$ given in Example \ref{ex:cantor}.}
            \label{fig:cantor}
        \end{figure}
        One could include an even broader class of convex sectors, e.g.~for negative
        $\alpha$-values, as transformations allow us to reduce these cases to those
        addressed directly by Theorem \ref{thm:main}. See the following discussion on
        future directions and Remark \ref{rem:omega}.
    \end{rem}
	
\section{Future Directions}
\label{sec:future}

    We now describe two open areas for future research.

    The Fueter and P{\'o}lya problem stated in Question \ref{question2},
    about the the existence of packing polynomials of higher degree on sectors, remains largely open.
    On $S(\infty)$ it is still open for degrees higher than 4,
    and it is indeed also an open problem for general sectors. 
    
    We now pose a question about general irrational sectors.

    \begin{question}
        Are there packing polynomials on sectors $S(\omega_1, \omega_2)$
        where $\omega_1$ and $\omega_2$ are both irrational?
    \end{question}

\begin{rem}
\label{rem:omega}
  	In this article we have focused on sectors of the type $S(\alpha)$
	rather than the general type given by two vectors $S(\omega_1, \omega_2)$.
	If one of $\omega_1$ or $\omega_2$ is rational, then the situations
	are equivalent, which we can see as follows.
	Assume that $\omega_1 = (r, s)$ is rational, by which we mean $r, s \in \N$, $\gcd(r,s) = 1$.
	We can then find integers $a, b \in \Z$ such that $a r + b s = 1$.
	The matrix $\begin{pmatrix} a & b\\ -s & r\end{pmatrix}$ has determinant 1
	and sends $(r,s)$ to $(1,0)$. If the transformation sends $\omega_2$ outside the first quadrant,
	we can apply transformations $\begin{pmatrix} 1 & 0\\ 0 & -1\end{pmatrix}$
	and $\begin{pmatrix} 1 & m\\ 0 & 1\end{pmatrix}$, $m \in \Z$,
	to accommodate that.
	So, if either of $\omega_1$ or $\omega_2$ is rational, the quadratic packing polynomials on $S(\omega_1, \omega_2)$ are known. If both $\omega_1$ and $\omega_2$ are irrational, then 
	it is unknown whether packing polynomials are possible.  
\end{rem}


\bibliographystyle{abbrv}
\bibliography{refs}

\end{document}